\documentclass[12pt, reqno]{amsart}
\usepackage{times}
\usepackage{amsmath,amsthm, amssymb, latexsym}
\usepackage{thmtools}
\usepackage{mathrsfs}
\newtheorem{theorem}{Theorem}[section]
\newtheorem{corollary}{Corollary}[section]
\newtheorem{lemma}{Lemma}[section]
\begin{document}
\title{Polynomial Criterion for Abelian Difference Sets }
\author{Pradipkumar H. Keskar and Priyanka Kumari}
\address{ Department of Mathematics, Birla Intitute of Technology and Science Pilani (Pilani Campus), Pilani 333031, India}
\email{keskar@pilani.bits-pilani.ac.in,\newline
priyanka.kumari@pilani.bits-pilani.ac.in}
\begin{abstract}
Difference sets are subsets of a group satisfying certain combinatorial property with respect to the group operation. They can be characterized using an equality in the group ring of the corresponding group. In this paper, we exploit the special structure of the group ring of an abelian group to establish a one-to one correspondence of the class of difference sets with specific parameters in that group with the set of all complex solutions of a specified system of polynomial equations. The correspondence also develops some  tests for a  Boolean function to be a bent function.


\end{abstract}

\maketitle

\section{Introduction}

For a finite group $\text{G}$ of order $v$ and nonnegative integers $k, \lambda$, a subset $\text{D}$ of $\text{G}$ is called a $(v, k, \lambda)$ \emph{difference set} in $\text{G}$ if for every $g \in \text{G} \setminus \{e\}$,
 $$|\{(d_1, d_2)\in \text{D} \times \text{D}\, :\, g=d_1d_2^{-1}\}|= \lambda\text{ and }|\text{D}|=k,$$
where $e$ is the identity of $\text{G}$. Moreover, if $\text{G}$ is abelian then $\text{D}$ is called an \emph{abelian difference set}.
 The notion of a difference set was introduced independently by J. Singer \cite{Sin} and R. C. Bose \cite{Bos} while investigating finite geometries and (statistical) design of experiments respectively. Later the ideas were found fruitful having several relations to areas such as coding theory and cryptography.  Bent functions, which are cryptographically significant, are characterized on page 95 of \cite{Sti} in terms of difference sets as follows . For an even positive integer $t > 2$, a Boolean function of $t$ variables \big(that is, a function from $\left(\mathbb{Z}/2\mathbb{Z}\right)^t$ to $\left(\mathbb{Z}/2\mathbb{Z}\right)$\big) is a bent function if and only if its support is a $\left(2^t, 2^{(t-1)} \pm 2^{(t-2)/2},2^{(t-2)} \pm 2^{(t-2)/2}\right)$ difference set in $\left(\mathbb{Z}/2\mathbb{Z}\right)^t$(where signs are chosen consistently).  Hence the construction, characterization, equivalence of difference sets is a useful exercise having applications for the analysis of bent functions.

 Among the different tools used to study difference sets are symmetric designs and group rings. In this paper, we will be concerned with the characterization of a difference set using an equality in a group ring, which we now describe.

 Let $\text{G}$ be a finite group and $\text{R}$ be a commutative ring with unit element $1$ different from its additive identity $0$. The \emph{group ring} $\text{RG}$ of $\text{G}$ over $\text{R}$ is the ring
 $$\text{RG} = \{ \sum_{g \in \text{G}}r_g g : r_g \in R\}$$
 where $\sum_{g \in \text{G}}r_g g = \sum_{g \in \text{G}}r_g^{*} g \iff r_g = r_g^{*}$ for all $g \in \text{G},$ with the addition defined by $\sum_{g \in \text{G}}r_g g + \sum_{g \in \text{G}}r_g^{*} g =\sum_{g \in \text{G}}(r_g+r_g^{*}) g$ and the multiplication defined by $\left(\sum_{g \in \text{G}}r_g g\right)  \left(\sum_{g \in \text{G}}r_g^{*} g\right) =\sum_{g \in \text{G}}(\sum_{xy=g}r_xr_y^{*}) g.$

For any $\text{D} \subset \text{G}$, we denote $\sum_{g \in \text{D}} g \in \text{RG}$ by $\text{D}$ again and $\sum_{g \in \text{D}} g^{-1} \in \text{RG}$ by $\text{D}^{(-1)}.$ The following is a characterization of a difference set in a finite group $\text{G}$.\\
 \underline{Group Ring Criterion :} Let $\text{G}$ be a finite group  of order $v$ and $ k, \lambda$ be nonnegative integers with $k \le v$. Then $\text{D} \subset \text{G}$   is a $(v, k, \lambda)$ difference set in $\text{G}$ if and only if as elements of $\mathbb{C}\text{G}$, we have $\text{D} \text{D}^{(-1)} = \lambda \text{G} + (k-\lambda)e$, where $e$ is the identity element of $\text{G}$.\\
(In most of the literature, for instance \cite{BJu}, the characterization is proved under extra assumption that $|\text{D}| = k$. This assumption can be seen to be superfluous, as it is implied by either of the above two conditions. It is enough to note that if $\text{D} \text{D}^{(-1)} = \lambda \text{G} + (k-\lambda)e$, then by comparing the coefficient of $e$ on both sides, we get $|\text{D}| = k$.)

In this paper, we exploit the structure of the group ring of an abelian finite group $\text{G}$ as an affine $\mathbb{C}$-algebra to obtain two algebraic criteria for a subset of $\text{G}$ to be a $(v, k, \lambda)$ difference set. The first criterion, Theorem 2.2 of Section 2, is in terms of an ideal membership problem. The second, Theorem 3.2 of Section 3, is via verification of some polynomial equations.  These processes can also be crystalised to give tests for a subset of an abelian group to be a $(v, k, \lambda)$ difference set.  The first test is verifiable using ideal theory or algebra softwares like Macaulay 2, while the second test is verifiable by explicit computations with complex numbers, especially the roots of unity.   Section 4 deals with generalization of these criteria to generalized difference sets. We illustrate the use of the criteria for difference sets through some examples in  Section 5.  The first two illustrate these tests for difference sets. The third illustrates how the polynomial criterion can be used to prove that a quadratic Boolean function is a bent function. In a future work we plan to explore the more of potential applications of these methods to the theory of bent functions.

Ideal theoretic methods, in particular Gr\"obner basis methods, are being widely applied to several problems in Science and Engineering, see \cite{Cox}. More specifically, we can find their applications to Combinatorics in \cite{Ron1}, \cite{Ron2}, \cite{Rob}, \cite{Ron3}. In this paper, we introduce these methods to study difference sets.

\section{Ideal membership problem for abelian difference sets }

Let $\text{G}$ be a finite abelian group. Then $\text{G} \cong C_{n_1} \times C_{n_2}\times \cdots \times C_{n_t}$ where $C_{n_l} = \left<g_l \right>$ is a cyclic group of order $n_l$. Let $ R = \mathbb{C}[X_1,\dots, X_t]$, $I$ be ideal $\left(X_1^{n_1}-1, \dots, X_t^{n_t}-1\right)$ of $R$ and $S = \{(i_1,\dots, i_t) \in \mathbb{Z}^t : 0 \le i_l \le n_l-1 \text{ for all } 1 \le l \le t\}$.  Regarding the structure of $\mathbb{C}\text{G}$, we have the following
\begin{theorem} Let $\text{G} \cong C_{n_1} \times C_{n_2}\times \cdots \times C_{n_t}$ where $C_{n_l}=\left<g_l\right> $ is a cyclic group of order $n_l$. Then the map  $\phi_\text{G}: \frac{R}{I} \to \mathbb{C}\text{G}$ defined by
$$\phi_\text{G}(f(X_1, \dots, X_t) +I)=f(g_1, \dots, g_t)$$
 is an isomorphism of $\frac{R}{I}$ onto $\mathbb{C}\text{G}$, where
\begin{eqnarray}
\nonumber f(g_1, \dots, g_t)&=&\sum_{(i_1,\dots,i_t)\in \mathbb{Z}^t}c_{i_1\cdots i_t}g_1^{i_1}\cdots g_t^{i_t} \\
 \nonumber \text{for }f(X_1,\dots,X_t)&=&\sum_{(i_1,\dots,i_t)\in \mathbb{Z}^t}c_{i_1\cdots i_t}X_1^{i_1}\cdots X_t^{i_t}\\
 \nonumber \text{ with }c_{i_1\cdots i_t} & \in & \mathbb{C}\text{ for all }(i_1,\dots,i_t)\in \mathbb{Z}^t.
    \end{eqnarray}
\end{theorem}

\noindent{\bf Proof :} First we show that $\phi_{\text{G}}$ is well defined. Assume $f_1(X_1,\dots,X_t)+I=f_2(X_1,\dots,X_t)+I$ for $f_1=f_1(X_1,\dots,X_t)$ and $f_2=f_2(X_1,\dots,X_t)\in R$.
Then $f_1-f_2\in I$ and hence $$f_1(X_1,\dots,X_t)-f_2(X_1,\dots,X_t)= \sum_{l=1}^t(X_l^{n_l}-1)u_l(X_1,\dots,X_t)$$ for $u_1(X_1,\dots,X_t),\dots,u_t(X_1,\dots,X_t)\in R$. This, in turn, implies that
$$f_1(g_1, \dots, g_t)-f_2(g_1,\dots,g_t)=\sum_{l=1}^t(g_l^{n_l}-1)u_l(g_1,\dots,g_t)=0$$ as $g_l^{n_l}-1=0$ for all $1 \le l \le t$.
Hence $\phi_\text{G}(f_1+I)=\phi_\text{G}(f_2+I)$, therefore $\phi_\text{G}$ is well defined.\\
Next, $\phi_\text{G}$ is clearly $\mathbb{C}-$algebra homomorphism onto $\mathbb{C}\text{G}$.\\
Now $\mathbb{C}\text{G}$ has dimension $n_1\cdots n_t$ as $\mathbb{C}-$vector space, as $\{g_1^{i_1}\cdots g_t^{i_t}: 0\le i_l\le n_l-1 \text{ for all }  1 \le l \le t\}$ is a basis for  $\mathbb{C}\text{G}$. Also $\frac{R}{I}$ has dimension $n_1\cdots n_t$ as $\{X_1^{i_1}\cdots X_t^{i_t}: 0\le i_l\le n_l-1 \text{ for all } 1 \le l \le t\}$ is a basis for $\frac{R}{I}$. This shows that $\phi_\text{G}$ is an isomorphism.   \qed
\\

As a consequence of Theorem 2.1, we can make several identifications. First, the group ring $\mathbb{C}\text{G}$ can be identified with the affine space $\mathbb{C}^{n_1 \cdots n_t}$ by identifying  $\sum_{(i_1,\dots,i_t)\in S}\alpha_{i_1\cdots i_t}g_1^{i_1}\cdots g_t^{i_t}$ with   $(\alpha_{i_1\cdots i_t} :(i_1,\dots,i_t)\in S)$ in a fixed order on $S$, say lexicographic order. Second, any subset $\text{T}$ of $\text{G}$ can be identified with the point in $\mathbb{C}^{n_1 \cdots n_t}$ corresponding to $\sum_{g \in \text{T}}g \in \mathbb{C}\text{G}$; it will be called {\em the point representation (or characteristic point) of } $\text{T}$ and it is a vertex of the unit hypercube of $\mathbb{C}^{n_1 \cdots n_t}$. Third, $\text{T}$ can also be represented by the unique polynomial $f=f(X_1,\dots, X_t) \in R$ such that $\phi_\text{G}(f+I)=\left(\sum_{g \in \text{T}}g\right)$ and either $f=0$ or $\text{deg}_{X_l}(f) < n_l$ for all $1 \le l \le t$. We call $f$ the {\em polynomial representation of }$\text{T}$ and denote it by $\rho_{\text{G}}(\text{T}) \text{ or } \rho_{\text{G}}(\text{T})(X_1,\dots ,X_t)$.

Here onwards, by using the isomorphism of $\text{G}$ with $C_{n_1} \times \dots \times C_{n_t}$ and fixed isomorphisms of $C_{n_l}$ with $\left(\frac{\mathbb{Z}}{n_l \mathbb{Z}}\right)$ for all $1 \le l \le t$, we will identify $\text{G}$ with $\prod_{l=1}^t \left(\frac{\mathbb{Z}}{n_l \mathbb{Z}}\right)$. Moreover for any $\text{T} \subset \text{G}$, we let $\text{T}^{*}=\{(i_1, \dots i_t) \in S: (i_1+n_1\mathbb{Z},\dots,i_t+n_t\mathbb{Z})\in \text{T}\}.$

The above relationships can be captured by the following equalities :
\begin{align*}
& \text{For any  T} \subset \text{G} \\
\rho_\text{G}(\text{T}) &= \sum_{(i_1, \dots, i_t) \in \text{T}^{*}}X_1^{i_1} \cdots X_t^{i_t} \tag{$2.1 *$}\\
\nonumber {} & = \sum_{(i_1, \dots, i_t) \in S} \alpha_{i_1 \cdots i_t} X_1^{i_1} \cdots X_t^{i_t}
\end{align*}
$ \text{where }\left(\alpha_{i_1\cdots i_t} : (i_1, \dots, i_t) \in S\right) \text{ is a point representation of T}$.

 Using these representations, Group Ring Criterion of Section 1 can be rephrased as an ideal membership problem (refer p. 94 of \cite{Cox}) in $\mathbb{C}\left[X_1, \dots, X_t\right]$ as follows.

\begin{theorem} Let $\kappa_{\text{D}} = \kappa_{\text{D}}(X_1, \dots, X_t) \in \mathbb{C}[X_1, \dots, X_t]$ be defined by
\begin{align*}
\kappa_{\text{D}} = &\left(\sum_{(i_1,\dots,i_t) \in \text{D}^{*}}X_1^{i_1}\cdots X_t^{i_t}\right)\left(\sum_{(i_1,\dots,i_t) \in \text{D}^{*}}X_1^{n_1-i_1}\cdots X_t^{n_t-i_t}\right) \\
&- \lambda \left(\sum_{(i_1, \dots, i_t) \in S}X_1^{i_1}\cdots X_t^{i_t}\right) - (k-\lambda).
\end{align*}
Then a subset $\text{D} \subset \text{G}$  is a $(v, k, \lambda)$ difference set in $\text{G}$ if and only if $\kappa_{\text{D}} \in I$.
 \end{theorem}

\noindent{\bf Proof :}
\begin{align*}
\phi_\text{G}(\kappa_\text{D}+I) =& \left(\sum_{(i_1,\dots,i_t) \in \text{D}^{*}}g_1^{i_1}\cdots g_t^{i_t}\right)\left(\sum_{(i_1,\dots,i_t) \in \text{D}^{*}}g_1^{n_1-i_1}\cdots g_t^{n_t-i_t}\right)\\
&-\lambda \left(\sum_{(i_1, \dots, i_t) \in S}g_1^{i_1}\cdots g_t^{i_t}\right) - (k-\lambda)\\
	=& \left(\sum_{(i_1,\dots,i_t) \in \text{D}^{*}}g_1^{i_1}\cdots g_t^{i_t}\right)\left(\sum_{(i_1,\dots,i_t) \in \text{D}^{*}}g_1^{-i_1}\cdots g_t^{-i_t}\right)\\
	&-\lambda \left(\sum_{(i_1, \dots, i_t) \in S}g_1^{i_1}\cdots g_t^{i_t}\right) - (k-\lambda)\\
	=& \text{DD}^{(-1)}-\lambda\text{G}-(k-\lambda).
	\end{align*}
Now
\begin{equation*}
\begin{split}
&\text{D is a }(v, k, \lambda) \text{ difference set in G} \\
 &\Leftrightarrow \text{DD}^{(-1)}-\lambda\text{G}-(k-\lambda)=0\\
 &\Leftrightarrow\phi_\text{G}(\kappa_\text{D}+I)=0 \Leftrightarrow \kappa_\text{D}+I=0\\  &\Leftrightarrow \kappa_\text{D}\in I.
 \end{split}
\end{equation*}

\qed

\noindent {\bf Note :} Alternatively we can write
\begin{align*}
\kappa_{\text{D}} = &\left(\sum_{(i_1,\dots,i_t) \in S}\alpha_{i_1\dots i_t}X_1^{i_1}\dots X_t^{i_t}\right)\left(\sum_{(i_1,\dots,i_t) \in S}\alpha_{i_1\dots i_t}X_1^{n_1-i_1}\dots X_t^{n_t-i_t}\right) \\
&- \lambda \left(\sum_{(i_1, \cdots, i_t) \in S}X_1^{i_1}\cdots X_t^{i_t}\right) - (k-\lambda).
\end{align*}
where $\alpha = (\alpha_{i_1\dots i_t} : (i_1,\dots,i_t) \in S)$ is the point representation of $\text{D}$.

\section{System of Polynomial equations for difference sets in an Abelian group}

The results of Section 2 can be rephrased to provide a criterion  for $(v, k, \lambda)$ difference sets in an abelian group of order $v$ in terms of some polynomial equations.  More specifically, we will find a set of polynomials in  $\mathbb{C}\left[\{A_{i_1 \cdots i_t} : (i_1, \dots, i_t) \in S\}\right]$ whose zero set in $\mathbb{C}^{n_1 \cdots n_t}$ is exactly the set of all the point representations of all $(v, k, \lambda)$ difference sets in $\text{G}$.

First some terminology and preparation. For any $J \subset \mathbb{C}\left[X_1, \dots, X_t\right]$, let $V(J) = \{(x_1, \dots, x_t)\in \mathbb{C}^t : f(x_1, \dots, x_t) =0 \text{ for all } f(X_1,\dots, X_t) \in J \}$. For any $W \subset \mathbb{C}^t$, let $I(W) = \{f(X_1, \dots, X_t) \in \mathbb{C}\left[X_1, \dots, X_t\right] : f(x_1, \dots, x_t) =0 \text{ for all } (x_1, \dots, x_t)\in W \}$. Then it can easily be seen that $I(W)$ is an ideal in $\mathbb{C}\left[X_1, \dots, X_t\right]$. For any ideal $J$ of $\mathbb{C}\left[X_1, \dots, X_t\right]$, the radical of $J$ is given by $\sqrt{J} = \{f=f(X_1, \dots, X_t) \in \mathbb{C}\left[X_1, \dots, X_t\right] : f^n \in J \text{ for some positive integer }n\}$. It can be seen that for any ideal $J$ of $\mathbb{C}\left[X_1, \dots, X_t\right]$, $\sqrt{J}$ is an ideal of $\mathbb{C}\left[X_1, \dots, X_t\right]$. An ideal $J$ of $\mathbb{C}\left[X_1, \dots, X_t\right]$ is called a radical ideal if $\sqrt{J}=J$. For the following famous theorem, see \cite{Cox}, p. 175.

\begin{theorem} {\bf (Hilbert Nullstellensatz)} If $J$ is any ideal of $\mathbb{C}\left[X_1, \dots, X_t\right]$ then $I(V(J))=\sqrt{J}$.
\end{theorem}

We will use the following corollary of Hilbert Nullstellensatz.
\begin{corollary} Let $f=f(X_1,\dots, X_t) \in \mathbb{C}\left[X_1, \dots, X_t\right]$ and $J$ be a radical ideal of $\mathbb{C}\left[X_1, \dots, X_t\right]$. Then $f \in J$ if and only if $f(x_1, \dots, x_t) = 0$ for all $(x_1, \dots, x_t) \in V(J)$.
\end{corollary}

In order to use Corollary 3.1, we prove the following
\begin{lemma} The ideal $I = \left( X_1^{n_1}-1, \dots, X_t^{n_t}-1\right)$ of $\mathbb{C}\left[X_1, \dots, X_t\right]$ is a radical ideal.
\end{lemma}

\noindent{\bf Proof :} Clearly $I \subset \sqrt{I}$. Now let $f(X_1, \dots, X_t) \in \sqrt{I}$ be any element. We want to show that $f=f(X_1, \dots, X_t) \in I$. We can write
$$f(X_1, \dots, X_t) = g(X_1, \dots, X_t) + r(X_1, \dots, X_t)$$
such that $g(X_1, \dots, X_t) \in I$ and $r(X_1, \dots, X_t) \in \mathbb{C}\left[X_1, \dots, X_t\right]$ satisfies $r(X_1, \dots, X_t)=0$ or $\text{deg}_{X_i}r(X_1, \dots, X_t) < n_i$ for all $i=1, 2, \dots, t$.

It is enough to show that $r(X_1, \dots, X_t)=0$. We can write
$$r(X_1,\dots,X_t) = \sum_{j=0}^{n_t-1}r_j(X_1, \dots, X_{t-1})X_t^j$$
 with $r_j(X_1, \dots, X_{t-1}) \in \mathbb{C}\left[X_1, \dots, X_{t-1}\right]$
 such that for any $0 \le j < n_t$, we have $r_j(X_1, \dots, X_{t-1})=0$ or $\text{deg}_{X_i}r_j(X_1, \dots, X_{t-1}) < n_i$ for all $1 \le i \le t-1$. Consider any $(\xi_1, \dots, \xi_{t-1}) \in \mathbb{C}^{t-1}$ with $\xi_i^{n_i}=1$ for all $1 \le i \le t-1$. Since $f \in \sqrt{I}=I(V(I))$ and for all $0 \le l \le n_t-1$, $\left(\xi_1, \dots, \xi_{t-1}, e^{\frac{2\pi i l}{n_t}}\right) \in V(I)$, we have $f\left(\xi_1, \dots, \xi_{t-1}, e^{\frac{2\pi i l}{n_t}}\right)=0$. Moreover, since $g \in I$,  $g\left(\xi_1, \dots, \xi_{t-1}, e^{\frac{2\pi i l}{n_t}}\right)=0$ for all $0 \le l \le n_t-1$. Thus $r\left(\xi_1, \dots, \xi_{t-1}, e^{\frac{2\pi i l}{n_t}}\right)=0$ for all $l =0, 1, \dots, n_t-1$. Assume $r(\xi_1, \dots, \xi_{t-1}, X_t) \ne 0$. Since $\text{deg}_{X_t}r(\xi_1, \dots, \xi_{t-1}, X_t) < n_t$ and it has $n_t$ distinct roots, we get a contradiction. Hence $r(\xi_1, \dots, \xi_{t-1}, X_t)=0$ and therefore for any $0 \le j < n_t$, we have $r_j(\xi_1, \dots, \xi_{t-1})=0$ for any $(\xi_1, \dots, \xi_{t-1}) \in \mathbb{C}^{t-1}$ with $\xi_i^{n_i}=1$ for all $i=1, 2, \dots, t-1$. Iterating the argument with $r$ replaced by $r_j$ etc, we get that $r=0$. \qed
\\

Now we are prepared to discuss the main result of this paper. Let $\Delta$ denote the polynomial ring $\mathbb{C}\left[X_1, \dots, X_t\right]\left[\{A_{i_1 \cdots i_t} :(i_1, \dots, i_t) \in S\}\right]$ in $n_1 \cdots n_t$ independent variables $A_{i_1 \cdots i_t} :(i_1, \dots, i_t) \in S$ over $\mathbb{C}\left[X_1, \dots, X_t\right]$
and let $U = \{(\xi_1,\dots,\xi_t) \in \mathbb{C}^t : \xi_i^{n_i}=1 \text{ for all }1 \le i \le t\}$.  To simplify the notation, let $A = \left(A_{i_1 \cdots i_t} : (i_1, \dots, i_t) \in S\right)$, $X = \left(X_1, \dots, X_t\right)$. Also if $\alpha_{i_1 \cdots i_t} \in \mathbb{C}$ for all $(i_1, \dots, i_t) \in S$, we let $\alpha = \left(\alpha_{i_1 \cdots i_t} : (i_1, \dots, i_t) \in S\right) \in \mathbb{C}^{n_1 \cdots n_t}$.
Then we have the following

\begin{theorem} Let $\Psi = \Psi(X, A) \in \Delta$ be defined by
\begin{align*}
\Psi  =& \left(\sum_{(i_1,\dots,i_t) \in S}A_{i_1\cdots i_t}X_1^{i_1}\cdots X_t^{i_t}\right)\left(\sum_{(i_1,\dots,i_t) \in S}A_{i_1\cdots i_t}X_1^{n_1-i_1}\cdots X_t^{n_t-i_t}\right) \\
&- \lambda \left(\sum_{(i_1, \dots, i_t) \in S}X_1^{i_1}\cdots X_t^{i_t}\right) - (k-\lambda)
\end{align*}
Then we have the following :
\begin{enumerate}
  \item{ For $\alpha = \left(\alpha_{i_1 \cdots i_t} : (i_1, \dots, i_t) \in S\right) \in \mathbb{C}^{n_1 \cdots n_t}$, $\alpha$ is a point representation of a subset of $\text{G}$ if and only if $\alpha$ satisfies the system $P_{i_1 \cdots i_t}(A) =0, (i_1, \dots, i_t) \in S$ of polynomial equations where for $(i_1, \dots, i_t) \in S, P_{i_1 \cdots i_t}(A)= A_{i_1 \cdots i_t}^2 - A_{i_1 \cdots i_t}$.}
  \item {For $\alpha = \left(\alpha_{i_1 \cdots i_t} : (i_1, \dots, i_t) \in S\right) \in \mathbb{C}^{n_1 \cdots n_t}$, $\alpha$ is a point representation of a $(v, k, \lambda)$ difference set in $\text{G}$} if and only if $\alpha$ satisfies the equations $P_{i_1 \cdots i_t}(A) =0\text{ for all } (i_1, \dots, i_t) \in S$, and $\Psi(\xi, A)=0$ for all $\xi = (\xi_1, \dots, \xi_t) \in U$.
\end{enumerate}
\end{theorem}

{\bf Proof :}Note that (1) follows, as for any $\alpha = \left(\alpha_{i_1 \cdots i_t} : (i_1, \dots, i_t) \in S\right) \in \mathbb{C}^{n_1 \cdots n_t}$,
\begin{equation*}
\begin{split}
& P_{i_1 \cdots i_t}(\alpha) =0 \text{ for all }(i_1, \dots, i_t) \in S \\
&\Leftrightarrow  \alpha_{i_1 \cdots i_t} \in \{0, 1\}\text{ for all }(i_1, \dots, i_t) \in S \\
& \Leftrightarrow  \alpha \text{ is a point representation of }\text{D} \subset \text{G} \text{ where } \\
&\text{D} = \{(i_1+n_1\mathbb{Z}, \dots, i_t+n_t\mathbb{Z}) : (i_1, \dots, i_t) \in S \text{ and }\alpha_{i_1 \cdots i_t}=1\}.
\end{split}
\end{equation*}

To prove (2), first note that $\kappa_{\text{D}}(X_1, \dots, X_t) = \Psi(X, \alpha)$ where $\alpha$ is the point representation of $\text{D} \subset \text{G}$. Thus, for any $\alpha \in \mathbb{C}^{n_1 \cdots n_t}$,

\begin{eqnarray}
\nonumber & \alpha \text{ is a point representation of a } (v, k, \lambda)\text{ difference set D in } \text{G}\\
\nonumber & \stackrel{\text{Theorem 2.2}}{\Longleftrightarrow} \alpha \text{ is a point representation of } \text{D} \subset \text{G} \text{ and }\\
\nonumber & \kappa_{\text{D}}(X_1, \dots, X_t) \in I \\
\nonumber & \Longleftrightarrow \alpha \text{ is a point representation of } \text{D} \subset \text{G} \text{ and }\\
\nonumber & \Psi(X, \alpha) \in I\\
\nonumber & \overset{\text{by Theorem 3.2 (1)}}{\Longleftrightarrow} P_{i_1 \cdots i_t}(\alpha) =0 \text{ for all }(i_1, \dots, i_t) \in S \text{ and }\\
\nonumber & \Psi(X, \alpha) \in I\\
\nonumber & \overset{\text{Theorem 3.1, Lemma 3.1}}{\Longleftrightarrow} P_{i_1 \cdots i_t}(\alpha) =0 \text{ for all }(i_1, \dots, i_t) \in S \text{ and }\\
\nonumber & \Psi(X, \alpha) \in I(V(I))\\
\nonumber & \stackrel{\text{since }U=V(I)}{\Longleftrightarrow} P_{i_1 \cdots i_t}(\alpha) =0 \text{ for all }(i_1, \dots, i_t) \in S \text{ and }\\
\nonumber & \Psi(\xi, \alpha)=0 \text{ for all } \xi \in U.
\end{eqnarray} \qed

\noindent {\bf Note :} The ideas of this section are analogous to interpolation of polynomials, in the sense that every postulation of a zero of a polynomial puts a condition on the parameters occurring in its coefficients.

\noindent {\bf Remark on $\mathbb{C}$-algebra homomorphisms :} For any $\xi = (\xi_1,\dots,\xi_t) \in U$, let $\theta_{(\xi_1, \dots, \xi_t)} : \mathbb{C}\left[X_1, \dots, X_t\right] \to \mathbb{C}$ be a $\mathbb{C}$-algebra homomorphism defined by $\theta_{(\xi_1, \dots, \xi_t)}\left(f(X_1, \dots, X_t)\right)=f(\xi_1, \dots, \xi_t)$ and let $\theta^{\Delta}_{(\xi_1, \dots, \xi_t)} : \Delta \to \Delta$ be the unique extension of $\theta_{(\xi_1, \dots, \xi_t)}$ to a $\mathbb{C}$-algebra homomorphism such that $\theta^{\Delta}_{(\xi_1, \dots, \xi_t)}(A_{i_1 \cdots i_t})=A_{i_1 \cdots i_t}$ for all $(i_1, \dots, i_t) \in S$. That is, for any $\Omega(X, A) \in \Delta$,
$$\theta^{\Delta}_{(\xi_1, \dots, \xi_t)}(\Omega(X, A) = \Omega(\xi, A).$$

Moreover, for any $\alpha \in \mathbb{C}^{n_1 \dots n_t}$, let $\tau_{\alpha} : \Delta \to \mathbb{C}\left[X_1, \dots, X_t\right]$ be a ring homomorphism defined by $\tau_{\alpha} \left(\Omega(X, A)\right) = \Omega(X, \alpha)$ for all $\Omega(X, A) \in \Delta$. Note that if $\alpha$ is the point representation of $\text{D} \subset \text{G}$ and $\xi = (\xi_1, \dots, \xi_t) \in U$,
\begin{eqnarray}
\nonumber  & \kappa_{\text{D}}(X_1, \dots, X_t)= \tau_{\alpha}(\Psi(X, A)),\\
\nonumber & \Psi(\xi, A) = \theta^{\Delta}_{(\xi_1, \dots, \xi_t)}(\Psi(X, A)),\\
\nonumber & \Psi(\xi, \alpha) = \theta_{(\xi_1, \dots, \xi_t)}\left(\tau_{\alpha}(\Psi(X, A))\right) = \theta_{(\xi_1, \dots, \xi_t)}\left(\kappa_{\text{D}}(X_1, \dots, X_t)\right)\\
\nonumber & \Psi(\xi, \alpha)=\tau_{\alpha} \left(\theta^{\Delta}_{(\xi_1, \dots, \xi_t)}(\Psi(X, A)\right) = \tau_{\alpha} \left(\Psi(\xi, A)\right).
\end{eqnarray}

Hence, in view of $(2.1*)$, it follows that if $\alpha$ is the point representation of $\text{D} \subset \text{G}$ and $\xi = (\xi_1, \dots, \xi_t) \in U$, then
\begin{align*}
\Psi(\xi, \alpha) = & \theta_{(\xi_1, \dots, \xi_t)}(\tau_{\alpha}(\Psi(X, A)))\\
= & \theta_{(\xi_1, \dots, \xi_t)}\left(\rho_{\text{G}}(\text{D})\right)\theta_{(\xi_1, \dots, \xi_t)}\left(\rho_{\text{G}}(\text{D}^{(-1)})\right)\\
& -  \lambda \ \theta_{(\xi_1, \dots, \xi_t)}\left(\rho_{\text{G}}(\text{G})\right) -(k-\lambda) \tag{$3.1 *$}.
\end{align*}

Alternatively, as a consequence of the conclusions of Theorem 3.2, for a subset $\text{D of G, D}$ is a $(v, k, \lambda)$ difference set if and only if
\begin{align*}
&\rho_{\text{G}}(\text{D})(\xi_1,\dots,\xi_t)\rho_{\text{G}}(\text{D}^{-1})(\xi_1,\dots,\xi_t)-\lambda \rho_{\text{G}}(\text{G})(\xi_1,\dots,\xi_t)-(k-\lambda)=0\\
& \text{for all } (\xi_1,\dots,\xi_t) \in U. \tag{$3.2 *$}
\end{align*}

\noindent {\bf Sharpening of Ryser Condition :} A necessary condition for the existence of $(v, k, \lambda)$ difference set is $\lambda(v-1)=k(k-1)$, as discovered by Ryser. Note that this condition is nothing but  $\theta^{\Delta}_{(\xi_1, \dots, \xi_t)}(\Psi)(\alpha)=0$ for $(\xi_1, \dots, \xi_t) =(1, \dots, 1)$, where $\alpha$ is the point representation of some set $\text{D} \subset \text{G}$ of size $k$. The condition is clearly not sufficient. However, when an abelian group $\text{G}$ of order $v$ is given as a direct sum of cyclic groups, the necessary as well as sufficient condition for existence of a $(v, k, \lambda)$ difference set in $\text{G}$ is the consistency of the system of equations  $P_{i_1 \cdots i_t}(A) =0\text{ for all } (i_1, \dots, i_t) \in S$ and $\theta^{\Delta}_{(\xi_1, \dots, \xi_t)}(\Psi)(A)=0$ for all $(\xi_1, \dots, \xi_t) \in U$ in the variables $A = \left(A_{i_1 \cdots i_t} : (i_1, \dots, i_t) \in S\right)$. In the parlance of Gr\"{o}bner bases (see page 171 of \cite{Cox}), the condition be reformulated as\\
\underline{Gr\"obner Basis Version of Existence Problem :} There exists a $(v, k, \lambda)$ difference set in $\text{G}$ if and only if the reduced Gr\"{o}bner basis (in any monomial order) of the ideal generated by $P_{i_1 \cdots i_t}(A) \text{ for all } (i_1, \dots, i_t) \in S$ and $\theta^{\Delta}_{(\xi_1, \dots, \xi_t)}(\Psi)(A)$ for all $(\xi_1, \dots, \xi_t) \in U$ in $\mathbb{C}\left[\{A_{i_1 \cdots i_t} : (i_1, \dots, i_t) \in S\}\right]$ is not equal to $\{1\}$.

\section{Generalization of the criteria}

The criteria developed in Sections 2 and 3 can be generalized to generalized difference sets in finite abelian groups. Following \cite{Cao}, we proceed to define a generalized difference set thus. For a finite group $\text{G}_0$ of order $v$, let $\text{D}_0, \text{M}_0 \subset \text{G}_0$ be such that $|\text{D}_0|=k >1, |\text{M}_0|>0$. For any $g \in \text{G}_0$, let $\lambda_g = |\{(d_1, d_2) \in \text{D}_0 \times \text{D}_0 : g = d_1d_2^{-1}\}|$. If $\lambda_1, \lambda_2$ are nonnegative integers, $\text{D}_0$ is called a $(v, |M_0|, k, \lambda_1, \lambda_2)$-{\em generalized difference set} of $\text{G}_0$ related to $\text{M}_0$ if for any nonidentity element $g \in \text{G}_0$
\begin{eqnarray}
\nonumber \lambda_g =\left\{
                       \begin{array}{ll}
                         \lambda_1, & \hbox{$\text{if }g \in \text{M}_0$;} \\
                         \lambda_2, & \hbox{$\text{if }g \not\in \text{M}_0$.}
                       \end{array}
                     \right.
\end{eqnarray}
This is a generalization of a difference set in the following sense. If $\text{M}_0=\{e\}$ for the identity element $e$ of $\text{G}_0$ and $\lambda_1 =0$ then any $(v, |M_0|, k, \lambda_1, \lambda_2)$-generalized difference set of $\text{G}_0$ related to $\text{M}_0$ is exactly a $(v, k, \lambda_2)$ difference set of $\text{G}_0$. Other special cases of generalized difference set give several important variations of a difference set in $\text{G}_0$. For instance, if $\text{M}_0$ is a subgroup of $\text{G}_0$ and $\lambda_1=0$, then $\text{D}_0$ is called a $(v, |M_0|, k, \lambda_2)$ \emph{relative difference set} of $\text{G}_0$ with relative to $\text{M}_0$. We say $\text{D}_0$ is a $(v, k, \lambda_1, \lambda_2)$ partial difference set in $\text{G}_0$ if $\text{M}_0=\text{D}_0$. Generalized difference sets are helpful in computation of autocorrelation of certain arrays, see \cite{Cao}. Partial difference sets have connections with strongly regular graphs, see \cite{Sma}. With this in mind, we state the polynomial criterion for the generalized difference sets as a consequence of the group ring criterion (\cite{Cao}, Theorem 2), the proof is analogous to  Theorem 3.2.

\begin{theorem}
Let $\text{M} \subset \text{G}$ and let $\Psi^{*}=\Psi^{*}(X,A) \in \Delta$ be defined by
\begin{eqnarray}
\nonumber  \Psi^{*}  =\left\{
                       \begin{array}{ll}
                        \left(\displaystyle \sum_{(i_1,\dots,i_t) \in S}A_{i_1\cdots i_t}X_1^{i_1}\cdots X_t^{i_t}\right)\left(\displaystyle \sum_{(i_1,\dots,i_t) \in S}A_{i_1\cdots i_t}X_1^{n_1-i_1}\cdots X_t^{n_t-i_t}\right) \\
- \lambda_1 \left(\displaystyle \sum_{(i_1, \dots, i_t) \in \text{M}^{*}} X_1^{i_1}\cdots X_t^{i_t}\right)\\-\lambda_2 \left(\displaystyle \sum_{(i_1, \dots, i_t) \in S \setminus \text{M}^{*}}X_1^{i_1}\cdots X_t^{i_t}\right) - (k-\lambda_1), \hspace{2.2cm} \textrm{if }e \in \text{M};\\
{}\\
\left(\displaystyle \sum_{(i_1,\dots,i_t) \in S}A_{i_1\cdots i_t}X_1^{i_1}\cdots X_t^{i_t}\right)\left(\displaystyle \sum_{(i_1,\dots,i_t) \in S}A_{i_1\cdots i_t}X_1^{n_1-i_1}\cdots X_t^{n_t-i_t}\right) \\
 -\lambda_1 \left(\displaystyle \sum_{(i_1, \dots, i_t) \in \text{M}^{*}}X_1^{i_1}\cdots X_t^{i_t}\right)\\-\lambda_2 \left(\displaystyle \sum_{(i_1, \dots, i_t) \in S \setminus \text{M}^{*}}X_1^{i_1}\cdots X_t^{i_t}\right) - (k-\lambda_2), \hspace{2.2cm} \textrm{if }e \not\in \text{M}. \\
                       \end{array}
                     \right.
\end{eqnarray}

  \noindent Then we have the following : \newline
  For $\alpha = \left(\alpha_{i_1 \dots i_t} : (i_1, \dots, i_t) \in S\right) \in \mathbb{C}^{n_1 \cdots n_t}$, $\alpha$ is a point representation of a $(v, |\text{M}|, k, \lambda_1, \lambda_2)$ generalized difference set in $\text{G}$ related to $\text{M}$ if and only if $P_{i_1 \cdots i_t}(\alpha) =0\text{ for all } (i_1, \dots, i_t) \in S$ and $\Psi^{*}(\xi,\alpha)=0$ for all $\xi=(\xi_1, \dots, \xi_t) \in U$.
\end{theorem}

\section{Illustrations of the criteria}

In this section, we illustrate the use of the criteria developed in sections 2 and 3. The purpose of the illustrations is just to show how the ideas developed in previous sections can potentially be used. As the outcomes of these illustrations are well known or can be proved by other means, some of the arguments, which are repetitive in nature, are left to the reader.
\\

\noindent{\bf Illustration  1 :} Let $\text{G}=\left(\frac{\mathbb{Z}}{4\mathbb{Z}}\right) \times \left(\frac{\mathbb{Z}}{4\mathbb{Z}}\right)$. We verify that the  subset $\text{D} = \{(0,1), (0,2), (0,3), (1, 0), (2, 0), (3,0)\}$ of $\text{G}$ is a $(16, 6, 2)$ difference set in $\text{G}$.

\noindent\underline{Method 1 : Ideal Membership Problem}

By Theorem 2.2, we need to see that $\kappa_{\text{D}}(X_1, X_2) \in I = (X_1^4-1, X_2^4-1)$. This can be done by verifying that the remainder of $\kappa_{\text{D}}(X_1, X_2)$ mod $I$ is $0$. The following Macaulay 2 code does this.

\begin{verbatim}
R= CC[x_1,x_2, MonomialOrder=> Lex]

I=ideal(x_1^4-1,x_2^4-1)

G=(x_1^3+x_1^2+x_1+1)*(x_2^3+x_2^2+x_2+1)

D=x_1^3+x_1^2+x_1+x_2^3+x_2^2+x_2

D1=(x_1^3*x_2^4+x_1^2*x_2^4+x_1*x_2^4+x_1^4*x_2^3
+x_1^4*x_2^2+x_1^4*x_2)

v=16

k=6

lambda=2

D*D1-(k-lambda)-lambda*G

oo%I
\end{verbatim}
\underline{Method 2 (Polynomial Criterion) :}

The point representation $\alpha=(\alpha_{ij} : (i,j) \in S)$ of $\text{D}$ is given by
\begin{eqnarray}
\nonumber \alpha_{01}= \alpha_{02}= \alpha_{03}=\alpha_{10}=\alpha_{20}=\alpha_{30} &=& 1 \\
\nonumber \alpha_{ij} &=& 0 \text{ elsewhere}.
\end{eqnarray}

As $\alpha_{ij} \in \{0, 1\}$, $P_{ij}(\alpha) =0, (i,j) \in S$.

Next, we verify that  $\theta^{\Delta}_{(\xi_1, \xi_2)}(\Psi)(\alpha)=0$ for all $(\xi_1, \xi_2) \in U$, where
 $U = \{\pm 1, \pm i\}^2$.

For sake of brevity, we verify one of the equations of Theorem 3.2 when $(\xi_1, \xi_2)=(i,-i) \in U$. The equations corresponding to the remaining 15 elements of $U$ can be verified to establish that $\text{D}$ is a  $(16, 6, 2)$ difference set in $\text{G}$.

Note that $\theta^{\Delta}_{(\xi_1, \xi_2)}(\Psi)(\alpha) =\Psi(i, -i,\alpha)$. Also
\begin{align*}
\Psi(u,v,\alpha) =& \left(u^0v^1+u^0v^2+u^0v^3+u^1v^0+u^2v^0+u^3v^0\right)\\
& \times \left(u^4v^3+u^4v^2+u^4v^1+u^3v^4+u^2v^4+u^1v^4\right)\\
&-2\left(\sum_{0 \le i \le 3, 0 \le j \le 3}u^iv^j\right) -(6-2).
\end{align*}

In particular, when $(u,v)=(i, -i)$, since
$$\sum_{0 \le r \le 3, 0 \le s \le 3}i^r(-i)^s = \left(\sum_{r=0}^3i^r\right)\left(\sum_{s=0}^3(-i)^s\right)=(0)(0),$$
we get
\begin{align*}
\Psi( i, -i, \alpha) =& (-i-1+i+i-1-i)(i-1-i+i-1-i)\\
&-2(0)-(6-2)\\ =& 0.
\end{align*}


The reader may also verify directly by definition that $\text{D}$ is a $(16, 6, 2)$ difference set in $\text{G}$.
\\

\noindent{\bf Illustration 2 :} Let $\text{D} = \{(0,1), (0,2), (0,3), (1, 0), (2, 0), (1, 1)\} \subset \text{G} =\left(\frac{\mathbb{Z}}{4\mathbb{Z}}\right) \times \left(\frac{\mathbb{Z}}{4\mathbb{Z}}\right)$. We can show that $\text{D}$ is not a $(16, 6, 2)$ difference set in $\text{G}$ just by verifying that $\Psi( 1, -1, \alpha) = -4 \ne 0$ for the point representation $\alpha$ of $\text{D}$.
\\

\noindent{\bf Illustration 3 (An application to bent functions)} This illustration provides the glimpse of the use of Theorem 3.2 to prove some results about bent functions. Let $t=2m$ for a positive integer $m$ and $\beta : \left(\mathbb{Z}/2\mathbb{Z}\right)^t \to \left(\mathbb{Z}/2\mathbb{Z}\right)$ be a Boolean function defined by $\beta(x_1, \dots, x_m, y_1, \dots, y_m) = \sum_{i=1}^m x_i y_i$ where $x_i, y_i \in \left(\mathbb{Z}/2\mathbb{Z}\right)$ are any elements for any $i = 1, \dots, m$. We illustrate the proof, using Theorem 3.2, that $\beta$ is a bent function with $|\text{D}| =2^{(t-1)}-2^{(t-2)/2}$ where $\text{D} = \text{support of }\beta$.

The proof can be given by induction on $m$. The cases $m=1, 2, 3$ can be verified easily. For the inductive step, let $m \geq 4$. Then we can write $m = m_1 + m_2$ with min$(m_1, m_2) \ge 2$. For any $i \in \{1,2\}$, let $t_i=2m_i$, $\text{G}_i = \left(\mathbb{Z}/2\mathbb{Z}\right)^{t_i}$ and let $\beta_i : \text{G}_i \to \left(\mathbb{Z}/2\mathbb{Z}\right)$ be defined by \newline $\beta_1(g_1)= \sum_{i=1}^{m_1} x_i y_i$ for any $g_1 \in \text{G}_1$, where $g_1 = (x_1, \dots, x_{m_1}, y_1, \dots, y_{m_1})$; $x_i, y_j \in \left(\mathbb{Z}/2\mathbb{Z}\right)$ for all $1 \leq i, j \leq m_1$ and $\beta_2(g_2)= \sum_{i=1}^{m_2} x_{m_1+i} y_{m_1+i}$ for any $g_2 = (x_{m_1+1}, \dots, x_{m}, y_{m_1+1}, \dots, y_{m}); x_i, y_j \in \left(\mathbb{Z}/2\mathbb{Z}\right)$ for all $m_1+1 \leq i, j \leq m$. By induction hypothesis, $\beta_1, \beta_2$ are bent functions with $|\text{D}_i| =2^{(t_i-1)}-2^{(t_i-2)/2}$ where $\text{D}_i = \text{support of }\beta_i$.

Let $\text{G}=\text{G}_1\times \text{G}_2$.  Identifying $g =(g_1, g_2) \in \text{G}$ with $(x_1, \dots, x_{m},y_1, \dots, y_{m}) \in \left(\mathbb{Z}/2\mathbb{Z}\right)^t$, $\text{G}$ gets identified with  $\left(\mathbb{Z}/2\mathbb{Z}\right)^t$.  Then $\beta(g_1, g_2) = \beta_1(g_1)+\beta_2(g_2)$ for any $g_1 \in \text{G}_1, g_2 \in \text{G}_2$ .
We need to show that $\beta$  is a bent function with domain $\text{G}$ and $|\text{D}| =2^{(t-1)}-2^{(t-2)/2}$ for $\text{D} = \text{support of }\beta$.

We start the proof with

\noindent\underline{Observations :}
\begin{enumerate}
	\item Let $\text{H}_i \subset \text{G}_i=\left(\mathbb{Z}/2\mathbb{Z}\right)^{r_i}$ for $i \in \{1, 2\}$ and let $\text{H} = H_1 \times H_2 \subset \text{G}_1 \times \text{G}_2$. Then
	
	$\rho_{\left(\text{G}_1 \times \text{G}_2\right)}\left(\text{H}_1 \times \text{H}_2\right) = \left(\rho_{\text{G}_1}(\text{H}_1)\right) \left(\rho_{\text{G}_2}(\text{H}_2)\right).$
	
	\item Let $\xi \in \mathbb{C}$ be a primitive $n^{\text{th}}$ root of unity. Then $\sum_{i=0}^{n-1} \xi^i = 0$.
	\item By Observations (1) and (2), if $(\xi_1, \dots, \xi_t) \in \{-1, 1\}^t \setminus \{(1, \dots, 1)\}$, then $\theta_{(\xi_1, \dots, \xi_t)}\left(\rho_{\text{G}}(\text{G})\right)=\sum_{(i_1, \cdots, i_t) \in S}\xi_1^{i_1}\cdots \xi_t^{i_t} =0$.
\end{enumerate}

Now $\text{D}_i$ is a $(2^{t_i}, 2^{(t_i-1)}-2^{(t_i-2)/2}, 2^{(t_i-2)}-2^{(t_i-2)/2})$ difference set for $i \in \{1, 2\}$. In order to show that $\beta$ is a bent function, it is enough to show that $\text{D}$ satisfies the conditions in Theorem 3.2 (b). Let $\alpha$ be the point representation of $\text{D}$. Then by Theorem 3.2 (a), $P_{i_1 \dots i_t}(\alpha) =0$ for all $(i_1, \dots, i_t) \in S$. Using $(3.1*)$, we will show that $\Psi(\xi, \alpha) =0$ for all $\xi=(\xi_1, \dots, \xi_t) \in U$.  Now $\text{D} = (\text{D}_1 \times \overline{\text{D}_2}) \cup (\overline{\text{D}_1} \times \text{D}_2)$, where $\overline{\text{D}_i}$ is the complement of $\text{D}_i$ in $\left(\mathbb{Z}/2\mathbb{Z}\right)^{t_i}$ for all $i \in \{1, 2\}$, and the union is disjoint. Hence by Observation (1)
\begin{align*}
\theta_{(\xi_1, \dots, \xi_t)}\left(\rho_{\text{G}}(\text{D})\right)& = \Big(\theta_{(\xi_1, \dots, \xi_{t_1})}\left(\rho_{\text{G}_1}(\text{D}_1)\right)\Big)\Big(\theta_{(\xi_{(t_1+1)}, \dots, \xi_{(t_1+t_2)})}\left(\rho_{\text{G}_2}(\overline{\text{D}_2})\right)\Big)\\
&+\Big(\theta_{(\xi_1, \dots, \xi_{t_1})}\left(\rho_{\text{G}_1}(\overline{\text{D}_1})\right)\Big)\Big(\theta_{(\xi_{(t_1+1)}, \dots, \xi_{(t_1+t_2)})}\left(\rho_{\text{G}_2}(\text{D}_2)\right)\Big).
\end{align*}
Therefore
\begin{align}
\nonumber{} &\Big( \theta_{(\xi_1, \dots, \xi_t)}\left(\rho_{\text{G}}(\text{D})\right)\Big)^2 \label{(*)}\\
\nonumber {} &= \Big(\theta_{(\xi_1, \dots, \xi_{t_1})}\left(\rho_{\text{G}_1}(\text{D}_1)\right)\Big)^2\Big(\theta_{(\xi_{(t_1+1)}, \dots, \xi_{(t_1+t_2)})}\left(\rho_{\text{G}_2}(\overline{\text{D}_2})\right)\Big)^2 \hfill{}\\
\nonumber {} &+\Big(\theta_{(\xi_1, \dots, \xi_{t_1})}\left(\rho_{\text{G}_1}(\overline{\text{D}_1})\right)\Big)^2\Big(\theta_{(\xi_{(t_1+1)}, \dots, \xi_{(t_1+t_2)})}\left(\rho_{\text{G}_2}(\text{D}_2)\right)\Big)^2 \\
\nonumber {} & + 2 \Big(\theta_{(\xi_1, \dots, \xi_{t_1})}\left(\rho_{\text{G}_1}(\text{D}_1)\right)\Big)\Big(\theta_{(\xi_{(t_1+1)}, \dots, \xi_{(t_1+t_2)})}\left(\rho_{\text{G}_2}(\overline{\text{D}_2})\right)\Big)\\
\nonumber {} & \Big(\theta_{(\xi_1, \dots, \xi_{t_1})}\left(\rho_{\text{G}_1}(\overline{\text{D}_1})\right)\Big)\Big(\theta_{(\xi_{(t_1+1)}, \dots, \xi_{(t_1+t_2)})}\left(\rho_{\text{G}_2}(\text{D}_2)\right)\Big) \tag{$5.1$}
\end{align}

We make three cases.

\noindent \underline{Case (i)} :
\begin{eqnarray*}
(\xi_1, \dots, \xi_{t_1}) \in \{-1, 1\}^{t_1} \setminus \{(1, \dots, 1)\}\text{ and }\\
(\xi_{(t_1+1)}, \dots, \xi_{(t_1+t_2)}) \in \{-1, 1\}^{t_2} \setminus \{(1, \dots, 1)\}.
\end{eqnarray*}
 In view of Observation (3), $\theta_{(\xi_1, \dots, \xi_{t_1})}(\rho_{\text{G}_1}(\overline{\text{D}_1})) = -\theta_{(\xi_1, \dots, \xi_{t_1})}(\rho_{\text{G}_1}(\text{D}_1))$ and $\theta_{(\xi_{(t_1+1)}, \dots, \xi_{(t_1+t_2)})}(\rho_{\text{G}_2}(\overline{\text{D}_2})) = -\theta_{(\xi_{(t_1+1)}, \dots, \xi_{(t_1+t_2)})}(\rho_{\text{G}_2}(\text{D}_2))$. Since $n_i =2$ for all $i =1, \dots, t_1$ in the equalities of Theorem 3.2, it follows that $\theta_{(\xi_1, \dots, \xi_{t_1})}(\rho_{\text{G}_1}(\text{D}_1^{(-1)})) = \theta_{(\xi_1, \dots, \xi_{t_1})}(\rho_{\text{G}_1}(\text{D}_1))$.

Similarly $\theta_{(\xi_{(t_1+1)}, \dots, \xi_{(t_1+t_2)})}(\rho_{\text{G}_2}(\text{D}_2^{(-1)}))= \theta_{(\xi_{(t_1+1)}, \dots, \xi_{(t_1+t_2)})}(\rho_{\text{G}_2}(\text{D}_2))$ and $\theta_{(\xi_1, \dots, \xi_t)}(\rho_{\text{G}}(\text{D}^{(-1)})) = \theta_{(\xi_1, \dots, \xi_t)}(\rho_{\text{G}}(\text{D}))$. Since each $\text{D}_i$ is a $(2^{t_i}, 2^{(t_i-1)}-2^{(t_i-2)/2}, 2^{(t_i-2)}-2^{(t_i-2)/2})$ difference set for $i \in \{1, 2\}$, as a consequence of $(3.1*)$ and Observation (3)
\begin{eqnarray*}
(\theta_{(\xi_1, \dots, \xi_{t_1})}(\rho_{\text{G}_1}(\text{D}_1)))^2 = 2^{(t_1-1)}-2^{(t_1-2)}\text{ and }\\  (\theta_{(\xi_{(t_1+1)}, \dots, \xi_{(t_1+t_2)})}(\rho_{\text{G}_2}(\text{D}_2)))^2 =2^{(t_2-1)}-2^{(t_2-2)}.
\end{eqnarray*}
Hence by $(5.1)$,
\begin{align*}
\Big( \theta_{(\xi_1, \dots, \xi_t)}\left(\rho_{\text{G}}(\text{D})\right)\Big)^2
=& 4\Big(\theta_{(\xi_1, \dots, \xi_{t_1})}\left(\rho_{\text{G}_1}(\text{D}_1)\right)\Big)^2
\left(\theta_{(\xi_{(t_1+1)}, \dots, \xi_{(t_1+t_2)})}\left(\rho_{\text{G}_2}(\text{D}_2)\right)\right)^2\\
=& 4 \left(2^{(t_1-1)}-2^{(t_1-2)}\right)\left(2^{(t_2-1)}-2^{(t_2-2)}\right)\\
=& \left(2^{(t-1)}-2^{(t-2)}\right) = k - \lambda
\end{align*}
where $k = 2^{t-1}-2^{(t-2)/2}, \lambda = 2^{t-2}-2^{(t-2)/2}$. Since $(\xi_1, \dots, \xi_t) \in \{-1, 1\}^t \setminus \{(1, \dots, 1)\}$, $\theta_{(\xi_1, \dots, \xi_t)}\left(\rho_{\text{G}}(\text{G})\right)=0$, and hence, by $(3.1*)$, $\Psi(\xi, \alpha)=0$ for the point representation $\alpha$ of $\text{D}$.
\\

\noindent \underline{Case (ii)} :
\begin{eqnarray*}
(\xi_1, \dots, \xi_{t_1})=(1, \dots, 1)\in \mathbb{C}^{t_1} \text{ and }\\
(\xi_{(t_1+1)}, \dots, \xi_{(t_1+t_2)}) \in \{-1, 1\}^{t_2} \setminus \{(1, \dots, 1)\}.
\end{eqnarray*}
 Then
\begin{align*}
\theta_{(\xi_1, \dots, \xi_{t_1})}\left(\rho_{\text{G}_1}(\text{G}_1)\right)=& 2^{t_1};\\
\theta_{(\xi_1, \dots, \xi_{t_1})}\left(\rho_{\text{G}_1}(\overline{\text{D}_1})\right)=&2^{t_1}-\theta_{(\xi_1, \dots, \xi_{t_1})}\left(\rho_{\text{G}_1}(\text{D}_1)\right);\\
\theta_{(\xi_{(t_1+1)}, \dots, \xi_{(t_1+t_2)})}\left(\rho_{\text{G}_2}(\text{G}_2)\right)=& 0;\\
\theta_{(\xi_{(t_1+1)}, \dots, \xi_{(t_1+t_2)})}\left(\rho_{\text{G}_2}(\overline{\text{D}_2})\right) =& -\theta_{(\xi_{(t_1+1)}, \dots, \xi_{(t_1+t_2)})}\left(\rho_{\text{G}_2}(\text{D}_2)\right).
\end{align*}

Hence by $(5.1)$,
\begin{align*}
&\Big( \theta_{(\xi_1, \dots, \xi_t)}\left(\rho_{\text{G}}(\text{D})\right)\Big)^2\\
&= \Big\{\theta_{(\xi_1, \dots, \xi_{t_1})}\left(\rho_{\text{G}_1}(\text{D}_1)\right)^2 + \left(2^{t_1}-\theta_{(\xi_1, \dots, \xi_{t_1})}\left(\rho_{\text{G}_1}(\text{D}_1)\right)\right)^2 \\
& -2 \theta_{(\xi_1, \dots, \xi_{t_1})}\left(\rho_{\text{G}_1}(\text{D}_1)\right)\left(2^{t_1}-\theta_{(\xi_1, \dots, \xi_{t_1})}\left(\rho_{\text{G}_1}(\text{D}_1)\right)\right)\Big\}\Big(\theta_{(\xi_{(t_1+1)}, \dots, \xi_{(t_1+t_2)})}\left(\rho_{\text{G}_2}(\text{D}_2)\right)\Big)^2.
\end{align*}
By (3.1*) and Observation (3), $\theta_{(\xi_{(t_1+1)}, \dots, \xi_{(t_1+t_2)})}\left(\rho_{\text{G}_2}(\text{D}_2)\right)^2 = 2^{t_2-1}-2^{t_2-2}$. Also $\theta_{(\xi_1, \dots, \xi_{t_1})}\left(\rho_{\text{G}_1}(\text{D}_1)\right) = |\text{D}_1| = 2^{t_1-1}-2^{(t_1-2)/2}$. By substituting these in the above expression and then expanding, we get
\begin{align*}
  & \left( \theta_{(\xi_1, \dots, \xi_t)}\left(\rho_{\text{G}}(\text{D})\right)\right)^2\\
 & =\left(2^{t_2-1}-2^{t_2-2}\right)\big[4\left(2^{t_1-1}-2^{(t_1-2)/2}\right)^2 +2^{2t_1} -4 (2^{t_1})\left(2^{t_1-1}-2^{(t_1-2)/2}\right)\big] \\
 &= \left(2^{t_2-1}-2^{t_2-2}\right)2^{t_1} = (k - \lambda)
\end{align*}
where $k = 2^{t-1}-2^{(t-2)/2}$ and $\lambda = 2^{t-2}-2^{(t-2)/2}$. Since $(\xi_1, \dots, \xi_t) \ne (1, \dots, 1)$, we see that $\theta_{(\xi_1, \dots, \xi_t)}\left(\rho_{\text{G}}(\text{G})\right)=0$ and hence by $(3.1*)$, $\Psi(\xi, \alpha)=0$ for the point representation $\alpha$ of $\text{D}$.

Similar argument works when $(\xi_1, \dots, \xi_{t_1})\in \{-1, 1\}^{t_2} \setminus \{(1, \dots, 1)\}$ and $(\xi_{(t_1+1)}, \dots, \xi_{(t_1+t_2)})=(1, \dots, 1) \in \mathbb{C}^{t_2} $.
\\

\noindent \underline{Case (iii)}:
\begin{eqnarray*}
(\xi_1, \dots, \xi_{t_1})=(1, \dots, 1) \in \mathbb{C}^{t_1} \text { and }\\
 (\xi_{(t_1+1)}, \dots, \xi_{(t_1+t_2)})=(1, \dots, 1) \in \mathbb{C}^{t_2}.
 \end{eqnarray*}
  Then we have
\begin{align*}
&\theta_{(\xi_1, \dots, \xi_{t_1})}\left(\rho_{\text{G}_1}(\text{D}_1)\right) = |D_1| = 2^{t_1-1}-2^{(t_1-2)/2}; \\
&\theta_{(\xi_{(t_1+1)}, \dots, \xi_{(t_1+t_2)})}\left(\rho_{\text{G}_2}(\text{D}_2)\right) = |D_2| = 2^{t_2-1}-2^{(t_2-2)/2};\\
&\theta_{(\xi_1, \dots, \xi_{t_1})}\left(\rho_{\text{G}_1}(\overline{\text{D}_1})\right) = |\overline{\text{D}_1}| = 2^{t_1}-2^{t_1-1}+2^{(t_1-2)/2}=2^{t_1-1}+2^{(t_1-2)/2}; \\
&\theta_{(\xi_{(t_1+1)}, \dots, \xi_{(t_1+t_2)})}\left(\rho_{\text{G}_2}(\overline{\text{D}_2})\right) = |\overline{\text{D}_2}| = 2^{t_2}-2^{t_2-1}+2^{(t_2-2)/2} = 2^{t_2-1}+2^{(t_2-2)/2}.
\end{align*}

Substituting in $(5.1)$, in view of $(3.1*)$, showing $\text{D}$ is a $(2^{t_1+t_2}, 2^{t_1+t_2-1}-2^{(t_1+t_2-2)/2}, 2^{t_1+t_2-2}-2^{(t_1+t_2-2)/2})$ difference set reduces to proving the following equality :
\\
If $t_1, t_2 \ge 4$ are even integers then
\begin{align}
\nonumber &\left(2^{t_2-1}+2^{(t_2-2)/2}\right)^2\left(2^{t_1-1}-2^{(t_1-2)/2}\right)^2 + \left(2^{t_1-1}+2^{(t_1-2)/2}\right)^2\left(2^{t_2-1}-2^{(t_2-2)/2}\right)^2\\
\nonumber &+2\left(2^{2t_1-2}-2^{t_1-2}\right)\left(2^{2t_2-2}-2^{t_2-2}\right)\\
\nonumber & = 2^{t_1+t_2}\left(2^{t_1+t_2-2}-2^{(t_1+t_2-2)/2}\right) + \left(2^{t_1+t_2-1}-2^{t_1+t_2-2}\right).
\end{align}

Now to prove this equality,
\begin{align*}
\nonumber L.H.S. = & \Big(2^{2t_1-2}+2^{t_1-2}-2^{\frac{3t_1}{2}-1}\Big)\Big(2^{2t_2-2}+2^{t_2-2}+2^{\frac{3t_2}{2}-1}\Big)\\
\nonumber & +\Big(2^{2t_1-2}+2^{t_1-2
}-2^{\frac{3t_1}{2}-1}+2\,2^{\frac{3t_1}{2}-1}\Big)\Big(2^{2t_2-2}+2^{t_2-2}-2^{\frac{3t_2}{2}-1}\Big)\\
\nonumber &+2\Big(2^{2t_1+2t_2-4}-2^{2t_1+t_2-4}-2^{t_1+2t_2-4}+2^{t_1+t_2-4}\Big).
\end{align*}

Simplifying further,
\begin{align*}
L.H.S. =& \Big(2^{2t_1-2} +2^{t_1-2}-2^{\frac{3t_1}{2}-1}\Big)\Big(2^{2t_2-1} +2^{t_2-1}\Big) + 2^{\frac{3t_1}{2}+2t_2-2}\\
&+ 2^{\frac{3t_1}{2}+t_2-2}
-2^{\frac{3t_1}{2}+\frac{3t_2}{2}-1}
+ 2^{2t_1+2t_2-3}
 -2^{2t_1+t_2-3}\\
 & -2^{t_1+2t_2-3} +2^{t_1+t_2-3}.
\end{align*}

More simplification gives
\begin{align*}
L.H.S. =& 2^{2t_1+2t_2-3} +2^{2t_1+t_2-3} + 2^{t_1+2t_2-3} + 2^{t_1+t_2-3} -2^{\frac{3t_1}{2}+2t_2-2}\\
& -2^{\frac{3t_1}{2}+t_2-2}
+2^{\frac{3t_1}{2}+2t_2-2} +2^{\frac{3t_1}{2}+t_2-2} - 2^{\frac{3t_1}{2}+\frac{3t_2}{2}-1}\\
&  + 2^{2t_1+2t_2-3} -2^{2t_1+t_2-3} -2^{t_1+2t_2-3} +2^{t_1+t_2-3}\\
= &2^{2t_1+2t_2-2} +2^{t_1+t_2-2} -2^{\frac{3t_1}{2}+\frac{3t_2}{2}-1}\\
=& R. H. S.
\end{align*}

This proves that $\beta$ is a bent function.                              \qed
\\

\noindent{\bf Conclusion :} In this paper, we have proved two algebraic criteria for a $(v, k, \lambda)$ difference set in a given abelian group of order $v$. Illustrations are provided indicating how they can be applied. Further applications are being planned.
\\

\noindent{\bf Acknowledgements :} Both the authors thank the support from FIST Programme vide SR/FST/MSI-090/2013 of DST, Govt. of India. The second author thanks UGC, Govt of India for the support under JRF Programme (SR. No. 2061540979, Ref. No. 21/06/2015(1)EU-V R. No. 426800). Both the authors thank S. Gangopadhyay and B. Mandal of IIT, Roorkee for stimulating discussions.


\begin{thebibliography}{9999}
\bibitem[1]{BJu}
Beth, T., Jungnickel, D., Lenz, H.
\emph{Design Theory, Volume 1}
Cambridge University Press, 1999.

\bibitem[2]{Bos}
Bose, R. C.
On the construction of balanced incomplete block designs,
Annals of Eugenics, vol. 9, p. 358-399, 1939 .

\bibitem[3]{Cao}
Cao, X., Sun, D.
Some nonexistence results on generalized difference sets,
Applied Mathematics Letters, vol. 21 no. 8, p. 797-802, 2008 .

\bibitem[4]{Cox}
Cox, D., Little, J. and O'Shea, D.
{\it Ideals,Varieties and Algorithms,}
 Springer Verlag, New York Inc, 1992 .

  \bibitem[5]{Ron1}
Felszeghy, B., R\'ath, B. and R\'onyai, L.
The lex game and some applications,
 Journal of Symbolic Computation, vol. 41, p. 663--681, 2006.

 \bibitem[6]{Ron2}
Felszeghy, B. and R\'onyai, L.
Some meeting points of Gr\"obner bases and combinatorics,
 {\it Algorithmic Algebraic Combinatorics and Gr\"obner bases (M.
Klin, G. A. Jones, A. Jurisic, M. Muzychuk, I. Ponomarenko editors)}, p. 207--227, Springer, 2009.

 \bibitem[7]{Rob}
 Kreuzner, M. and Robbiano, L.
 {\it Computational Commutative Algebra 2}, Springer, 2005.

\bibitem[8]{Sma}
Ma, S.
A survey of partial difference sets,
 Designs, Codes and Cryptography, vol. 4, p. 221-261, 1994.

  \bibitem[9]{Ron3}
 R\'onyai, L. and M\'esz\'aros, T.
Some combinatorial applications of  of Gr\"obner bases,
 {\it Algebraic Informatics (Franz Winkler ed.), 4th International Conference, CAI
2011, Linz, Proceedings; Springer LNCS 6742}, p. 65--83, Springer-Verlag 2011.

\bibitem[10]{Sin}
Singer, J.
A theorem in finite projective geometry and some applications to number theory,
Transactions of American Mathematical Society, vol. 43, p. 377-385, 1938

\bibitem[11]{Sti}
Stinson, D. R.
{\it Combinatorial Designs : Construction and Analysis,}
 Springer Verlag, New York Inc, 2004 .

\end{thebibliography}
\end{document}